\documentclass[final]{article}

\usepackage{amsmath,amssymb}
\usepackage{amsfonts}
\usepackage{eufrak}
\usepackage{amsthm}
\usepackage{amsbsy}
\usepackage{latexsym}

\newcommand{\R}{\mathbb{R}}

\usepackage[utf8]{inputenc}
\usepackage[english]{babel}
\usepackage{url}
\usepackage{hyperref}
\usepackage{graphicx}
\usepackage{listings}
\usepackage{courier}
\usepackage{color}
\usepackage{amsmath}
\usepackage{amssymb}
\usepackage[T1]{fontenc}

\usepackage{tikz,pgfplots}

\pgfdeclarelayer{background}
\pgfdeclarelayer{outerneuron}
\pgfsetlayers{background,outerneuron,main}

\begin{document}

\pagestyle{headings}

\title{Hybrid solver methods for ODEs: Machine-Learning combined with standard methods.}

\author{
J\"urgen Geiser
\thanks{Data-Science, Freelancer, Berlin, E-mail: juergen.geiser@geiserdatascience.com}
}

\maketitle

\begin{abstract}
In this article, we consider combined standard and machine learning methods to solve ODEs and PDEs.
We deal with the minimisation problems for the machine learning algorithms and standard discretization methods, which are related
to Runge-Kutta methods and finite difference methods.
We show, that we could solve the ODEs with additional ML methods, e.g., feedforward network, such that it will accelerate the
solver process.

\end{abstract}



\bigskip

\section{Introduction}

We are motivated to add ML concepts into our solver methods for ODEs and PDEs.
While, we have a known part of the differential equations, e.g., linear part,
we could apply fast classical numerical or analytical solvers, the other more
delicate or less known parts, e g., nonlinear part, stochastic parts or less structured parts,
while we have training sets, should be done with neutral networks.
Here, standard solvers will fail, while the new neutral network solvers could find the
solver structures.

The paper is organised as follows. Mathematical models of the wave equations are introduced in Section \ref{s:model}. 
In section \ref{ml}, we give an overview to the proposed machine learning and numerical analysis of the combined methods.
First ideas are presented and finally, we discuss our future works in the area of ML problems in section \ref{summary}.

\section{Mathematical model}
\label{s:model}

We deal with general differental equation, e.g, dynamical system or semi-discretised PDE, which is given in the following manners:

\begin{itemize}
\item Linear and nonlinear ODEs or PDEs:
\begin{eqnarray}
\label{dyn_1}
&& \frac{\partial {\bf y}}{\partial t} = A u + {\bf F}({\bf y}(t), {\bf u}(t)) , \; t \in [0, T] , \\
&& {\bf y}(0) = {\bf y}_0, 
\end{eqnarray}
where ${\bf y}$ are the states variable and ${\bf F}$ is a nonlinear vectorial function, ${\bf v}_0$ is a constant initial vector and ${\bf u}$ are external variables.
Further the linear part is given with the matrix-operator $A \in \R^m \times \R^m$.

\item Deterministic and stochastic ODEs or PDEs:
\begin{eqnarray}
\label{dyn_1_1}
&& d {\bf y} = A u dt  + {\bf F}({\bf y}(t), {\bf u}(t)) dW , \; t \in [0, T] , \\
&& {\bf y}(0) = {\bf y}_0, 
\end{eqnarray}
where ${\bf y}$ are the states variable and ${\bf F}$ is a nonlinear vectorial function, ${\bf v}_0$ is a constant initial vector and ${\bf u}$ are external variables.
The stochastic part is given with the variable $W$ with $W(t)$ is  a Wiener-process with $\Delta W = W(t^{n+1}) - W(t^n) = \sqrt{\Delta t} \xi$ and $\xi$ is a normal
distributed variable with $\xi \approx N(0, 1)$. 
Further the linear part is given with the matrix-operator $A in \R^m \times \R^m$.

\end{itemize}

We apply ML and standard schemes, while we apply the fast standard methods to the
linear or deterministic parts of the differential equations and the 
nonlinear or stochastics and even more complex parts of the differential equation to the ML processes.

For the training of the dynamical system (\ref{dyn_1_1}) with the combination of the
classical solvers, we assume a numerical solution given with finite differences as:

\begin{itemize}
\item Linear-Nonlinear Splitting:
\begin{eqnarray}
\label{dyn_2}
&& \hspace{-1cm} {\bf y}_{num}(t + \Delta t) = \exp(A \Delta t) {\bf y}_{num}(t) + \Delta t \; \exp(A \Delta t) \; {\bf F}({\bf y}_{num}(t), {\bf u}(t)) , \\
&& {\bf y}_{num}(0) = {\bf y}_0, 
\end{eqnarray}
where $\Delta t$ is the time-step (sufficient small, e.g., CFL condition of an explicit Euler method)
and ${\bf y}_{linear}(t  + \Delta t) = \exp(A \Delta t) {\bf y}_{num}(t) $.

\item Linear-Nonlinear Splitting of higher order (Strang-Splitting approach)

\end{itemize}

\begin{figure}[ht]
\begin{center}  
\includegraphics[width=10.0cm,angle=-0]{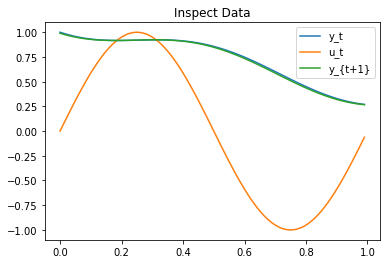}
\end{center}
\caption{\label{result_1} Results of the forward Euler method.}
\end{figure}

As a result, you will obtain the output, see Figure \ref{result_1}.

\section{Machine-Learning methods}
\label{ml}

We apply a feed-forward network with one hidden layer, which can be done with the
software-package {\it Tensorflow}).

\begin{figure}[ht]
\begin{center}  
\includegraphics[width=5.0cm,angle=-0]{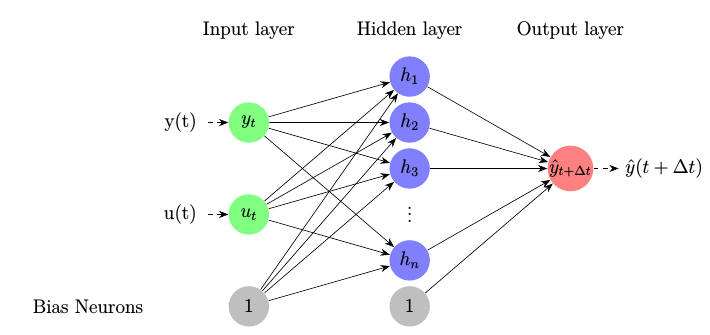}
\end{center}
\caption{\label{feedforward} Sketch of the feedforward network to approximate the linear system.}
\end{figure}

We like to simulate {\it nonlinear dynamical systems} with a feedforward
network, which is given in the Figure \ref{feedforward}.

\begin{figure}[ht]
\begin{center}  
\includegraphics[width=5.0cm,angle=-0]{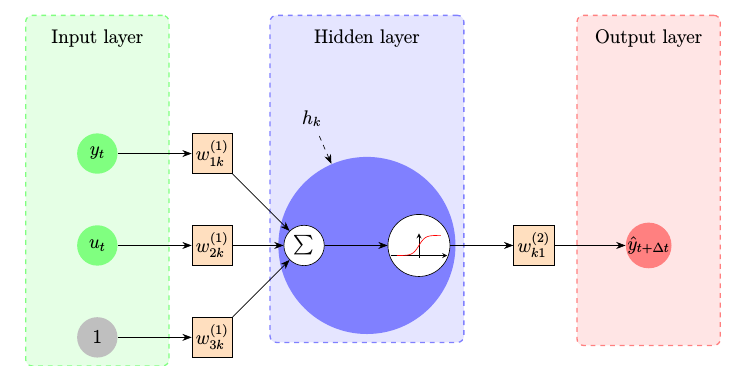}
\end{center}
\caption{\label{feedforward-hidden} Sketch of a hidden neuron $k$ of the feedforward network.}
\end{figure}

Taking a closer look on the dataflow (Figure \ref{feedforward-hidden}), we can see that the input of the k-th neuron in a hidden layer can be computed by taking the scalar product of a weight vector and the output of the previous layers neurons.
The neuron's output is simply the application of its activation function onto this scalar product.
We will use the sigmoid function as an activation function, which is given by
\begin{eqnarray}
\label{dyn_1}
\sigma(x) = \frac{1}{1+\exp({-x})} 
\end{eqnarray}

Looking at all neurons in the l-th layer at once this process can be described with a matrix-vector multiplication between a weighting matrix $w^{(l)}$ and the output of the previous layer's neurons.

We are interested in finding an optimal weighting of our neurons for a given measurement of goodness.
For simplicity we take the \textit{mean squared error}.
Therefore we have formulated the training process of a neural network as a minimization problem:
\begin{eqnarray}
\label{dyn_3}
&& min || {\bf \hat{y}} -  {\bf y}_{num} \|^2
\end{eqnarray}
where $|| \cdot ||$ is the vector-norm, ${\bf \hat{y}}$ is the result of the
neutral network and ${\bf y}_{num}$ is the numerical solution.

In general the training examples are based on analytical or numerical solutions
of such nonlinear/linear dynamical systems.

We can write down the network as an analytical function:
\begin{eqnarray}
\label{dyn_4}
&& {\bf \hat{y}}(t_k) = w^{(2)}(\sigma(w^{(1)} ({\bf y}(t_{k-1}), {\bf u}(t_{k-1}), 1)^T), 1)^T  , 
\end{eqnarray}
where $t_k = \Delta t \; k $ and $k = 1, \ldots, K$ with $T = K \; \Delta t$,
and \\
$({\bf \hat{y}}(t_{k-1}), {\bf \hat{u}}(t_{k-1}), 1)= (\hat{y}_1(t_{k-1}), \ldots, \hat{y}_m(t_{k-1}), \hat{u}_1(t_{k-1}) \ldots, \hat{u}_m(t_{k-1}), 1)$

Note that in this formulation we apply $\sigma$ to each element individually.
$w^{(1)} \in \mathbb{R}^{3 \times n}$ and $w^{(2)}\in \mathbb{R}^{(n+1) \times 1}$ are the matrices containing the neural network's weights between the layers with the matrix 
\begin{eqnarray}
A = \left( \begin{array}{c c c c}
A_{1,1} & A_{1,2} & \ldots & A_{1, 2m} \\
A_{2,1} & A_{2,2} & \ldots & A_{2, 2m} \\
\vdots & \vdots & \ddots & \vdots \\
A_{m,1} & A_{m,2} & \ldots & A_{m, 2m} 
 \end{array}\right) .
\end{eqnarray}

In the following, we have the task to implement the neutral network in tensorflow and to train it with the samples generated in task 1.
To accomplish this we need several tensorflow functions. For simplicity here is a list of all relevant functions.
\begin{itemize}
\item Placeholder
\item matmul
\item Variable
\item zeros
\item add
\item nn.sigmoid
\item reduce\_mean
\item reduce\_sum
\item squared\_difference
\item train.GradientDescentOptimizer
\item global\_variables\_initializer
\item Session
\end{itemize}

Start by skimming through the documentation of these functions.

\begin{figure}[ht]
\begin{center}  
\includegraphics[width=10.0cm,angle=-0]{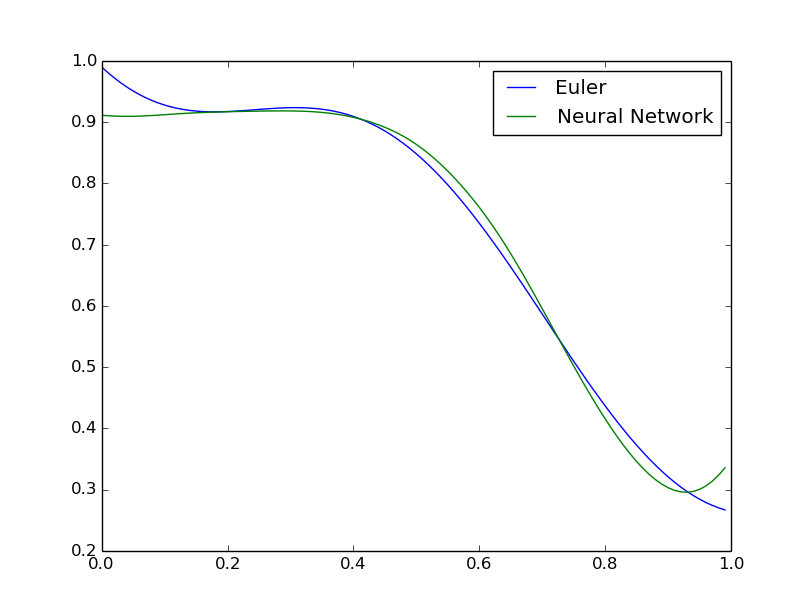}
\end{center}
\caption{\label{result_2} Results of the Feed forward network with one hidden-layer.}
\end{figure}

Initialize the bias weights with zeros and all other weights with normal distributed values.
We apply 10 hidden neurons, which are already enough for this task.
As a result, you will obtain the output, see Figure \ref{result_2}.

We could apply the softwaren and go on with the numerical approach of the
ODE-system with:
\begin{eqnarray}
\label{dyn_6}
&& {\bf y}_{num}(t + \Delta t) =  {\bf y}_{num}(t) + \Delta t \; {\bf F}({\bf y}_{num}(t), {\bf u}(t)) , \\
&& {\bf y}_{num}(0) = {\bf y}_0, 
\end{eqnarray}
where $\Delta t$ is the time-step (sufficient small, e.g., CFL condition of an explicit Euler method).

We have the following examples:
\begin{itemize}
\item Linear first decay with external force:
\begin{eqnarray}
\label{dyn_2}
&& \frac{\partial y}{\partial t} =  \lambda y + u(t) , t \in [0, T]\\
&& y(0) = y_0, 
\end{eqnarray}
where $T = 1.0$.

Numerical solution is given as:
\begin{eqnarray}
\label{dyn_21}
&& y_{num}(t + \Delta t) =  y_{num}(t) + \Delta t \; \lambda y_{num}(t) + u(t) , \\
&& y_{num}(0) = 1.0, 
\end{eqnarray}
where $\lambda = - 0.1$ and $u(t) = \sin(2 \pi t)$ and CFL-condition $\Delta t \le \frac{1}{\lambda}$.
\end{itemize}

The activation function is given as
\begin{eqnarray}
\label{dyn_1}
\sigma(x) = \frac{1}{1+\exp({-x})} 
\end{eqnarray}
where $y_{a,j}$ is the output of the summation function.

We use a first layer FF-NN, which is given as Figure \ref{feedforward}.

For training we again minimize the mean square error by applying an evolution-strategy. 
We start by instanciating $N$ differently parameterised neural networks, called \textit{population}.
The different neural networks acre called individuals.
As a fitness function we use the mean squared error, which we aim to minimize.
Then we repeat the following procedure $M$ times, where M is a large integer:
\begin{enumerate}
    \item Evaluate the fitness function of each individual
    \item Sort the population according to their fitness
    \item For each individual (indexed by $i$ trough the previous sorting) in the complement of the fittest 1/5th of the population, take the parametrization individual indexed by $i$ modulo 5, add gaussian noise to it and store this as the new parametrization of $i$. The gaussian noise must be sufficiently small, such that the reparametrization is ''somewhat smooth''.
\end{enumerate}
Taking the parametrization of the fittest individual after sufficient iterations ($M \approx 100$ for $N = 250$ and multiplying the values of the gaussian distribution by a factor of $0.05$) usually yields a satisfying parametrization.

As a result, you will obtain the same output as in the tensorflow implementation.

\begin{figure}[ht]
\begin{center}
\includegraphics[width=10.0cm,angle=-0]{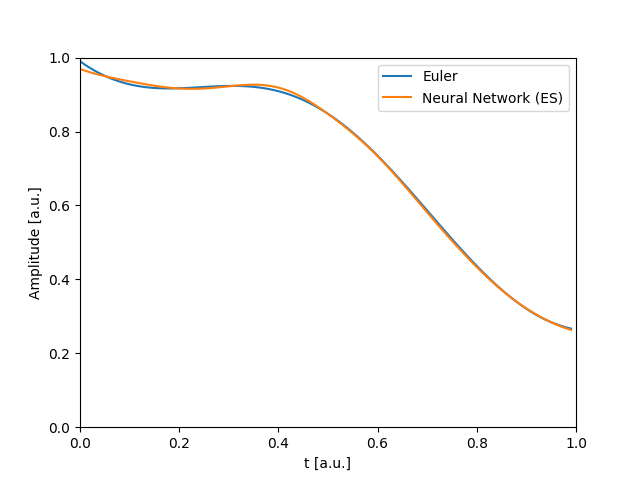}
\end{center}
\caption{\label{result_3} Results of the neural network with $M = 100$ and $N = 250$ compared to the Euler solution.}
\end{figure}

A comparison of the results from the neral network with the best fitness function and the numerical result is shown in Figure \ref{result_3}.
\subsection{Derivation of the minimisation problem}

We apply the derivation of the gradient of the MSE.

As the previous example underlines so-called blackbox optimization methods usually take a lot of time.
Since we have a ''nice'' representation of our optimization problem we can use gradient information to speed-up finding optima, as the gradient gives us information about the local slope of each dimension and therefore about the direction towards a minimum.
Therefore take the gradient of the MSE w.r.t. the weights.

The training of the network is done via the stochastic descent gradient.

We apply the gradient and train the network based on the MSE with stochastic descent gradient. 
In stochastic gradient descent we simply take the gradient of the MSE over a single function evaluation:

$$
w \leftarrow w - \eta \nabla \| {\bf \hat{y}}(t) -  {\bf y}_{num}(t) \|^2
$$

which is applied repeatedly for all samples.

As a result, you will again obtain the same output as in the tensorflow implementation.

\subsection{Optimization problem non-convex}

 Given a locally minimal parametrization of the neural network. 
We can obtain more optimal solutions by permutating this parametrization properly.
For example, if can take in this example two neurons from the hidden layer and find another optimal solution by swapping all of their weights.
When taking into account that there can be $n$ neurons in the hidden layer we receive $2^n - 1$ equivalent optimal parametrizations.
I will skip the proof for existance of at least one local minimum, as it should be intuitively clear in the case of the mean squared error.
Since we can construct (for the case $n>1$) multiple local minima, the problem cannot be convex.

This is usually not a problem, as these local minima are equivalently good.
Additionally to this, global minima can be obstructive, as we usually train neural networks to generalize to unseen data and not to overfit the samples.
Therefore even saddle points can be sufficiently good, see \cite{dl2016}

We have a convex problem, based on the well-defined linear equation system and the convex solution-space.
While this is true for the solutions of the linear dynamical system, this does not hold for neural networks.

\subsection{Derive the gradients for a general feedforward network}

Let us now gerneralize by deriving the gradient of feedforward neural networks with arbitrary (sufficiently differentiable) activation functions and an arbitrary amount of layers.
I.e. how do we obtain and explicit formula for
$$
\frac{\partial MSE}{\partial w^{(k)}_{ij}}
$$
w.r.t. to the adjacent layers.

For more details of the solution, see \cite{dl2016}.
Such ideas are found in in the deep learning book, see \cite{dl2016}.

\subsection{Validation of the network}

Here, we validate the network for more or less unseen data, while the standard solvers, e.g., in our case the
forward Euler method will fail.

We assume to have some more training dates, e.g., the analytical solution of the underlying 
equation:
\begin{eqnarray}
y(t) = \exp(\lambda t) y(0) ,
\end{eqnarray}
where $\lambda$ is our underlying parameter.

Now, we could train the NN with more accurate results, while the standard classical solver will have a problem with 
to coarse time-step.

In the following figure, we see the benefit of the {\it more} trained NN.

\subsection{Application of the ML solver for the PDE}

We deal with the PDE, which is given as:
\begin{align}
& \frac{\partial y}{\partial t}  = D \frac{\partial^2 y}{\partial x^2} , (x, t) \in [0, 10]\times [0, 1] , \\
& y(x, 0) = y_0(x) , \;  x \in [0, 10] , \; \mbox{initial condition} , \\
& y(0, t) = y(1, t) = 0 , \; t \in [0, 1] , \; \mbox{boundary condition} ,
\end{align}
where $D = 0.1$ and $u_0(x) = \left\{ \begin{array}{l l} 1 , & x \in [4.5, 5.5] , \\ 0 , & \mbox{else} \end{array} \right.$

Based on the spatial discretisation, we obtain a system of ODEs:
\begin{eqnarray}
\label{dyn_1_2}
&& \frac{\partial {\bf y}}{\partial t} = A {\bf y}, \; t \in [0, T] , \\
&& {\bf y}(0) = {\bf y}_0, 
\end{eqnarray}
where we have ${\bf y}(t) = (y_1, \ldots, y_I)^t \in \R^I$. Further, we have:
\begin{align}
A &= \frac{D}{\Delta x^2} \begin{pmatrix}
-2 & 1 & 0 & 0 & \cdots \\ 
1 & -2 & 1 & 0 &  \\ 
0 & 1 & -2 & 1 &  \\
\vdots & & & & \ddots
\end{pmatrix} \in \R^{I \times I}.
\end{align}

Further, we have the spatial step $\Delta x = 0.1$, means $I = 100$ and the time-step $\Delta t = 0.05$, means $N = 20$.

For the training of the dynamical system (\ref{dyn_1_2}), we assume a numerical solution given with finite differences as:
\begin{eqnarray}
\label{dyn_2}
&& {\bf y}_{num}(t + \Delta t) =  {\bf y}_{num}(t) + \Delta t \; A {\bf y}_{num}(t) , \\
&& {\bf y}_{num}(0) = {\bf y}_0, 
\end{eqnarray}
where $\Delta t \le \frac{\Delta x^2}{2 D}$ is the time-step (sufficient small, e.g., CFL condition of an explicit Euler method).

We apply to minimize the solution of the forward-network solver:
\begin{eqnarray}
\label{dyn_1}
&& || {\bf \hat{y}}(t) -  {\bf y}_{num}(t) \| \le min , \; t \in [0, T] ,
\end{eqnarray}
where $|| \cdot ||$ is the vector-norm, ${\bf \hat{y}}$ is the result of the
neutral network and ${\bf y}_{num}$ is the numerical solution.

We can write down the network as an analytical function:
\begin{eqnarray}
\label{dyn_1}
&& {\bf \hat{y}}(t_k) = w^{(2)}(\sigma(w^{(1)} ({\bf y}(t_{k-1}), 1)^T), 1)^T  , 
\end{eqnarray}
where $t_k = \Delta t \; k $ and $k = 1, \ldots, K$ with $T = K \; \Delta t$,
and \\
$({\bf \hat{y}}(t_{k-1}), 1)= (\hat{y}_1(t_{k-1}), \ldots, \hat{y}_m(t_{k-1}), 1)$

Note that in this formulation we apply $\sigma$ to each element individually.
$w^{(1)} \in \mathbb{R}^{2 \times N}$ and $w^{(2)}\in \mathbb{R}^{(N+1) \times 1}$ are the matrices containing the neural network's weights between the layers.

For the improved training, we assume to have the analytical solution:
\begin{eqnarray}
\label{dyn_2}
&& {\bf y}_{ana}(t) =  \exp(A t) {\bf y}_0, 
\end{eqnarray}
while $A$ is the spatial discretised matrix.

For such an improve training, we validate the NN and obtain more accurate results as for the standard solvers.

\section{Conclusions and Discussions}
\label{summary}

We have proposed a mixture between standard and ML algorithms, that helps to solve ODE and PDE problems with respect to new
methods related to ML algorithms.
Such combination can help to accelerate the solver processes and solve standard ODEs with the help of
ML methods. In future, we analyse more complex application of such problems.

\bibliographystyle{plain}

\begin{thebibliography}{10}



\bibitem{gei_07}
J.~Geiser,
\newblock Iterative Operator-Splitting Methods with higher order Time-Integration Methods and Applications for Parabolic Partial Differential
Equations,
\newblock J. Comput. Appl. Math., accepted, June 2007.

\bibitem{rack2020}
C.~Rackauckas et. al.
\newblock{\em Universal Differential Equations for Scientific Machine Learning}.
\newblock Arxiv - Jan 2020.


\bibitem{dl2016}
I.~Goodfellow, Y.~Bengio and A.~Courville,
\newblock{\em Deep Leaerning}.
\newblock MIT Press, 2016.
  
\end{thebibliography}

\end{document}